\documentclass[10pt]{amsart}
\usepackage{mathptmx}
\usepackage{amsmath}
\usepackage{amssymb}
\usepackage{array}
\usepackage{geometry}
\usepackage[bookmarks=true,colorlinks=true, pdfstartview=FitV, linkcolor=black, citecolor=blue, urlcolor=black]{hyperref}

\usepackage{color}
\definecolor{DarkRed}{rgb}{0.55,.00,0.2}
\definecolor{DarkGrey}{rgb}{0.35,.35,0.35}

\theoremstyle{definition}

\theoremstyle{remark}

\numberwithin{equation}{section}

\begin{document}

\title{ New index transforms  of the Lebedev- Skalskaya type}

\author{S. Yakubovich}
\address{Department of Mathematics, Faculty of Sciences,  University of Porto,  Campo Alegre str.,  687; 4169-007 Porto,  Portugal}
\email{ syakubov@fc.up.pt}


\date{\today}
\maketitle

\markboth{\rm \centerline{ S.  Yakubovich}}{}
\markright{\rm \centerline{Lebedev- Skalskaya Type Transforms }}

\begin{abstract}   New index transforms, involving the real part of the modified Bessel  function of the first  kind as the kernel  are  considered.   Mapping properties such as the boundedness and  invertibility are investigated for these  operators  in the Lebesgue  spaces.  Inversion theorems are proved.    As an interesting application,  a solution of the initial   value problem for the second order partial differential equation, involving the Laplacian,  is obtained.  It is noted,  that the corresponding operators  with the imaginary part of the modified Bessel function of the first kind lead to the familiar Kontorovich- Lebedev transform and its inverse. 
\end{abstract}

\section{Introduction and preliminary results}

Let $f(x),  \ g(\tau),\  x \in \mathbb{R}_+,\  \tau \in  \mathbb{R}$ be  complex -valued functions.  The main goal of this paper is to investigate mapping properties of the following index transforms of the Lebedev-Skalskaya type  \cite{yak}, involving  the modified Bessel function of the first   kind  in  the kernel,  namely,
$$(Ff) (\tau) = {\sqrt\pi\over  \cosh(\pi\tau)} \int_0^\infty e^{-x/2} {\rm Re} \left[ I_{i\tau} \left({x\over 2}\right) \right]  f(x)dx,   \quad    \tau \in \mathbb{R}, \eqno(1.1)$$
$$(G g) (x) =  \sqrt\pi \  e^{-x/2} \int_{-\infty}^\infty   {\rm Re}  \left[ I_{i\tau} \left({x\over 2}\right)\right] \   { g(\tau) d\tau\over \cosh(\pi\tau)} ,   \quad   x  \in \mathbb{R}_+, \eqno(1.2)$$
where $i$ is the imaginary unit and ${\rm Re}$ denotes the real part of a complex -valued function.  The modified Bessel function of the first kind  $I_\nu(z)$ \cite{erd}, Vol. II  satisfies the differential equation
$$  z^2{d^2u\over dz^2}  + z{du\over dz} - (z^2+\nu^2)u = 0.\eqno(1.3)$$
It has the asymptotic behavior 
$$ I_\nu(z) = {e^{z} \over \sqrt{2\pi z} }   [1+ O(1/z)], \ z \to \infty,\  - {\pi\over 2} < \arg z < {3\pi\over 2} \eqno(1.4)$$
and 
$$I_\nu(z) = O( |z|^{{\rm Re} \nu} ), \ z \to 0.\eqno(1.5)$$
The modified Bessel function of the first kind has the following series representation  
$$I_{\nu}(z)= \sum_{k=0}^\infty {(z/2)^{2k+\nu} \over k! \Gamma(k+\nu+1)},\   z,  \nu \in \mathbb{C},\eqno(1.6)$$
where $\Gamma(z)$ is Euler's gamma function \cite{erd}, Vol. I.  Hence with the reduction  formula for the gamma function \cite{erd}, Vol. I we find for ${\rm Re} \nu \ge 0$
$$| \Gamma(k+\nu+1) |= | \Gamma(\nu+1) (1+\nu)(2+\nu)\dots (k+\nu)| \ge k!  | \Gamma(\nu+1)|$$
and we derive from (1.6)
$$|I_{\nu}(z)| \le   e^{-  {\rm Im} \nu \arg z } \sum_{k=0}^\infty {(|z|/2)^{2k+  {\rm Re} \nu} \over k! |\Gamma(k+\nu+1)|} \le 
e^{-  {\rm Im} \nu \arg z } \  { \left(|z| /2\right)^{  {\rm Re} \nu} \over  | \Gamma(\nu+1)| } \sum_{k=0}^\infty {(|z|/2)^{2k} \over (k!)^2}$$
$$\le \  e^{|z| -  {\rm Im} \nu \arg z } \  { \left(|z| /2\right)^{  {\rm Re} \nu} \over  | \Gamma(\nu+1)| } ,$$
coming up to the following inequality for the modified Bessel function of the first kind 
$$|I_{\nu}(z)| \le \  \left({|z|\over 2}\right)^{  {\rm Re} \nu} \  \frac {e^{|z| -  {\rm Im} \nu \arg z }}{ | \Gamma(\nu+1)| },\ z,  \nu \in \mathbb{C} .\eqno(1.7)$$
In the mean time, another solution of the equation (1.3) is the Macdonald function $K_\nu(z)$ \cite{erd}, Vol. II
$$K_\nu(z) = {\pi\over 2\sin (\pi\nu)} \left[ I_{-\nu}(z)- I_\nu(z)\right].\eqno(1.8)$$
In particular, letting $\nu=i\tau,\ \tau \in \mathbb{R},  z=x >0$, we find from (1.8)
$$K_{i\tau}(x)= - {\pi\over \sinh(\pi\tau) } {\rm Im} \left[ I_{i\tau} \left(x\right) \right]\eqno(1.9)$$
and this function is the kernel of the classical Kontorovich-Lebedev transform \cite{yak}, \cite{yal}.   Correspondingly, taking into account the value
$$|\Gamma(1+i\tau)|=  \sqrt {{\pi\tau\over \sinh(\pi\tau)}},$$
inequality (1.7) takes the form
$$ |I_{i\tau}(x)| \le \  e^x\  \sqrt{{\sinh(\pi\tau)\over \pi\tau}},\   x >0,  \tau  \in \mathbb{R}.\eqno(1.10)$$

On the other hand, appealing to relation (8.4.22.5) in \cite{prud}, Vol. III, we find the following Mellin-Barnes representation for the modified Bessel function of the first kind 
$$ e^{-x/2}  I_{\nu} \left({x\over 2}\right) = {1\over 2\pi\sqrt\pi  i} \int_{\gamma-i\infty}^{\gamma +i\infty} \frac {\Gamma(s+\nu)\Gamma(1/2-s)}{\Gamma (\nu+1-s)} x^{-s} ds,\ - {\rm Re} \nu < \gamma < {1\over 2}.\eqno(1.11)$$

{\bf Lemma 1}. {\it Let $x >0,  \tau \in \mathbb{R}$. Then the following Mellin-Barnes integral representation of the kernel in $(1.1), (1.2)$ takes place}
$${\sqrt\pi\over  \cosh(\pi\tau)} \  e^{-x/2} {\rm Re} \left[ I_{i\tau} \left({x\over 2}\right) \right] \equiv  \Psi_\tau(x) =
{1\over 2\pi  i} \int_{\gamma-i\infty}^{\gamma +i\infty} \frac {\Gamma(s+ i\tau)\Gamma(s-i\tau) \Gamma(1/2-s)}{\Gamma(s)\Gamma (1-s)} x^{-s} ds,\   0 < \gamma < {1\over 2}.\eqno(1.12)$$
\begin{proof} In fact, taking (1.11) with $\nu=i\tau$, we have  
$$ e^{-x/2}  {\rm Re} \left[ I_{i\tau} \left({x\over 2}\right)\right] = {1\over 4\pi\sqrt\pi  i} \int_{\gamma-i\infty}^{\gamma +i\infty}  \frac {\Gamma(s+ i\tau)\Gamma (1-i\tau-s) +  \Gamma(s- i\tau)\Gamma (i\tau+1-s)}{ \Gamma(1-i\tau-s)\Gamma (i\tau+1-s)}$$
$$\times \ \Gamma(1/2-s)  x^{-s} ds.\eqno(1.13)$$
Meanwhile, employing the reflection formula for the gamma function \cite{erd}, Vol. I and elementary trigonometric formulae, we find 
$${1\over 4\pi\sqrt\pi  i} \int_{\gamma-i\infty}^{\gamma +i\infty}  \frac {\Gamma(s+ i\tau)\Gamma (1-i\tau-s) +  \Gamma(s- i\tau)\Gamma (i\tau+1-s)}{ \Gamma(1-i\tau-s)\Gamma (i\tau+1-s)}\ \Gamma(1/2-s)  x^{-s} ds$$
$$= {\sqrt \pi \over 4\pi  i} \int_{\gamma-i\infty}^{\gamma +i\infty}  \frac {[\sin(\pi(s+i\tau)) +  \sin(\pi(s-i\tau))]\  \Gamma(1/2-s) }{\sin(\pi(s+i\tau) )\sin(\pi(s-i\tau)) \Gamma(1-i\tau-s)\Gamma (i\tau+1-s)}\  x^{-s} ds$$

$$= {\cosh(\pi\tau) \over 2\pi^2\sqrt \pi  i} \int_{\gamma-i\infty}^{\gamma +i\infty}   \Gamma(s+i\tau)\Gamma (s- i\tau)\sin(\pi s) \Gamma(1/2-s) \  x^{-s} ds$$

$$=  {\cosh(\pi\tau) \over 2\pi \sqrt \pi  i} \int_{\gamma-i\infty}^{\gamma +i\infty} \frac {\Gamma(s+ i\tau)\Gamma(s-i\tau) \Gamma(1/2-s)}{\Gamma(s)\Gamma (1-s)} x^{-s} ds.$$
Hence, combining with (1.13), we  arrive at the equality (1.12), completing the proof of Lemma 1. 
\end{proof}

Equality (1.12) can be used to calculate for all $x >0$ the Fourier cosine transform \cite{tit} by $\tau$ of the kernel in (1.1).

{\bf Lemma 2}. {\it Let $x, y  >0$. Then the following index integral converges absolutely and has the value
$$ \int_0^\infty {\cos(\tau y) \over  \cosh(\pi\tau)} \  {\rm Re} \left[ I_{i\tau} \left({x\over 2}\right) \right] d\tau =  \sqrt {{x\over \pi}} \  e^{x/2}  \cosh(y/2) \   {}_1F_1 \left( 1; \  {3\over 2}; \   -x \cosh^2 \left({y \over 2}\right) \right),\eqno(1.14)$$
where ${}_1F_1 \left( a; \  b;  \  z \right) $ is the Kummer confluent hypergeometric function \cite{erd}, Vol. I.}

\begin{proof}    The absolute convergence of the integral in the left-hand side of (1.14) follows immediately from the inequality (1.10). Hence, multiplying both sides of (1.12) by $\cos(\tau y)$ and integrating over $\mathbb{R}_+$, we appeal to reciprocal formulae via the Fourier cosine transform (cf. formula (1.104) in \cite{yak}) 
$$\int_0^\infty  \Gamma\left(s + i\tau\right)  \Gamma\left(s - i\tau \right)  \cos( \tau y) d\tau
= {\pi\over 2^{2s}}  {\Gamma(2s) \over \cosh^{2s}(y/2)},\ {\rm Re}\ s > 0,\eqno(1.15)$$
$$  \Gamma\left(s + i\tau \right)  \Gamma\left(s - i\tau\right)  
=   { \Gamma(2s)  \over 2^{2s-1}}  \int_0^\infty   {\cos(\tau y)  \over \cosh^{2s} (y/2)} \ dy,\eqno(1.16)$$ 
and reverse the order of integration in the obtained right-hand side of (1.12).  It is indeed possible due to the Fubini theorem and the inequalities 
$$\left| \Gamma\left(s + i\tau\right)  \Gamma\left(s - i\tau\right)  \right|
\le {|\Gamma(2s+1)|\over \tau^2} \left[ c_1+ c_2 |s|\  \right],\  {\rm Re}\ s > 0,\  \tau \in  \mathbb{R} \backslash \{0\} ,\eqno(1.17)$$
$$ \left| \Gamma\left(s + i\tau\right)  \Gamma\left(s - i\tau\right)  \right| \le |\Gamma (2s)|
B(\gamma,\gamma), \eqno(1.18)$$
where $c_1, c_2$ are absolute positive constants and  $B(a,b)$ is Euler's beta-function \cite{erd}, Vol. I.  We note that inequality (1.17) can be easily obtained via integration  by parts  twice in (1.16). Hence with the use of  the Stirling asymptotic formula for the gamma-function \cite{erd}, Vol. I    it guarantees the absolute convergence of the corresponding iterated integral. Thus,  recalling (1.15), the duplication formula for the gamma- function \cite{erd}, Vol. I, we calculate the integral via Slater's theorem \cite{prud}, Vol. III  in terms of the Kummer function to find
$$ \int_0^\infty {\cos(\tau y) \over  \cosh(\pi\tau)} \  {\rm Re} \left[ I_{i\tau} \left({x\over 2}\right) \right] d\tau = e^{x/2}
{\sqrt \pi \over 2\pi   i} \int_{\gamma-i\infty}^{\gamma +i\infty} \frac {\Gamma(2s) \Gamma(1/2-s)}{\Gamma(s)\Gamma (1-s)} \left( 4x \cosh^2(y/2) \right)^{-s} ds $$
$$=   {e^{x/2}   \over 4 \pi   i} \int_{\gamma-i\infty}^{\gamma +i\infty} \frac {\Gamma(s+ 1/2) \Gamma(1/2-s)}{\Gamma (1-s)} \left( x \cosh^2(y/2) \right)^{-s} ds $$
$$=  \sqrt {{x\over \pi}} \  e^{x/2}  \cosh(y/2)   {}_1F_1 \left( 1; {3\over 2};  -x \cosh^2\left({y\over 2}\right) \right).$$

\end{proof} 

{\bf Corollary 1}.  {\it Let $x >0, \tau  \in \mathbb{R}.$  The kernel $(1.12)$ of the Lebedev-Skalskaya type transform $(1.1)$ has the integral representation 
$$  \Psi_\tau(x)  =    {2\sqrt x\over \pi}\ \int_0^\infty \cos(\tau y) \cosh(y/2) \   {}_1F_1 \left( 1; \  {3\over 2}; \   -x \cosh^2 \left({y \over 2}\right) \right) dy,\eqno(1.19)$$
where the integral converges absolutely. }

\begin{proof}  The proof is immediate from the inversion formula for the Fourier cosine transform and asymptotic behavior of the Kummer function at infinity (cf. \cite{erd}, Vol. I)
$${}_1F_1 \left( a; \   b; \   - z \right)  = O( z^{-a}),\ z \to +\infty,$$
which guarantees the absolute convergence of the integral (1.19) because the integrand is continuously differentiable as a function of $y \in \mathbb{R}_+$.
\end{proof} 

Employing the Mellin-Barnes representation (1.12) of the kernel  $\Psi_\tau(x)$, we will derive an ordinary differential equation whose particular solution is  $\Psi_\tau(x)$. Precisely, it is given by 

{\bf Lemma 3}.    {\it The kernel $\Psi_{\tau}(x)$ is a fundamental solution of the following second  order differential equation with variable  coefficients}
$$x^2 {d^2 \Psi_{\tau} \over dx^2} + x (1+x)  {d \Psi_{\tau} \over dx} + \left({x\over 2} +\tau^2\right) \Psi_{\tau}  = 0. \eqno(1.20)$$

\begin{proof}  Recalling the Stirling asymptotic formula for the gamma-function, we see that for each $\tau \in \mathbb{R}$ the integrand in (1.12) behaves as ($s= \gamma +it$)
$$ \frac {\Gamma(s+ i\tau)\Gamma(s-i\tau) \Gamma(1/2-s)}{\Gamma(s)\Gamma (1-s)} = 
e^{-\pi |t|/2} \  |t|^{\gamma - 1},\quad  |t| \to \infty.\eqno(1.21)$$
This circumstance   allows to differentiate repeatedly with respect to $x$ under the integral sign in (1.12).   Hence with the reduction formula for the gamma-function, a simple change of variables and the Cauchy residue theorem, we obtain the chain of equalities
$$x {d\over dx} x {d\over dx}\  \Psi_{\tau} =  {1\over 2\pi  i} \int_{\gamma-i\infty}^{\gamma +i\infty} \frac {s^2\  \Gamma(s+ i\tau)\Gamma(s-i\tau) \Gamma(1/2-s)}{\Gamma(s)\Gamma (1-s)} x^{-s} ds $$
$$= - \tau^2  \Psi_{\tau}(x)  + {1\over 2\pi  i} \int_{\gamma-i\infty}^{\gamma +i\infty} \frac {\Gamma(1+ s+ i\tau)\Gamma(1+ s-i\tau) \Gamma(1/2-s)}{\Gamma(s)\Gamma (1-s)} x^{-s} ds $$
$$= -  \tau^2  \Psi_{\tau}(x) + {1\over 2\pi  i} \int_{1+ \gamma-i\infty}^{1+ \gamma +i\infty} \frac {\Gamma( s+ i\tau)\Gamma( s-i\tau) \Gamma(3/2-s)}{\Gamma(s-1)\Gamma (2-s)} x^{1-s} ds   $$
$$=  -  \tau^2  \Psi_{\tau}(x)  - {1\over 2\pi  i} \int_{1+ \gamma-i\infty}^{1+ \gamma +i\infty} \frac {\Gamma( s+ i\tau)\Gamma( s-i\tau) \Gamma(1/2-s)(1/2-s) }{\Gamma(s)\Gamma (1-s) } x^{1-s} ds  $$
$$ =   \left( {3 x \over 2}-  \tau^2 \right) \   \Psi_{\tau}(x) - {1\over 2\pi  i} \int_{ \gamma-i\infty}^{\gamma +i\infty} \frac {\Gamma( s+ i\tau)\Gamma( s-i\tau) \Gamma(1/2-s)(2-s) }{\Gamma(s)\Gamma (1-s) } x^{1-s} ds  $$
$$ =   \left( {3 x \over 2}-  \tau^2 \right) \   \Psi_{\tau}(x) -  {d\over dx}\ \left( x^2\   \Psi_{\tau} (x) \right).$$
Hence  after fulfilling the differentiation  we arrive at (1.20), completing the proof of Lemma 3.

\end{proof}

\section {Boundedness  and inversion properties of the index transform (1.1)}

In order to investigate the mapping properties of the index transform (1.1) we will use the Mellin transform technique developed in \cite{yal}.   Precisely, the Mellin transform is defined, for instance, in  $L_{\nu, p}(\mathbb{R}_+),\ 1 < p \le 2$ (see details in \cite{tit}) by the integral  
$$f^*(s)= \int_0^\infty f(x) x^{s-1} dx,\eqno(2.1)$$
 being convergent  in mean with respect to the norm in $L_q(\nu- i\infty, \nu + i\infty),\   q=p/(p-1)$.   Moreover, the  Parseval equality holds for $f \in L_{\nu, p}(\mathbb{R}_+),\  g \in L_{1-\nu, q}(\mathbb{R}_+)$
$$\int_0^\infty f(x) g(x) dx= {1\over 2\pi i} \int_{\nu- i\infty}^{\nu+i\infty} f^*(s) g^*(1-s) ds.\eqno(2.2)$$
The inverse Mellin transform is given accordingly
 $$f(x)= {1\over 2\pi i}  \int_{\nu- i\infty}^{\nu+i\infty} f^*(s)  x^{-s} ds,\eqno(2.3)$$
where the integral converges in mean with respect to the norm  in   $L_{\nu, p}(\mathbb{R}_+)$
$$||f||_{\nu,p} = \left( \int_0^\infty  |f(x)|^p x^{\nu p-1} dx\right)^{1/p}.\eqno(2.4)$$
In particular, letting $\nu= 1/p$ we get the usual space $L_1(\mathbb{R}_+)$.  Further, denoting by $C (\mathbb{R}_+)$ the space of bounded continuous functions, we have 

{\bf Theorem 1.}   {\it The index transform  $(1.1)$  is well-defined as a  bounded operator $F: L_1(\mathbb{R}_+) \to C (\mathbb{R})$ and the following norm inequality takes place
$$||Ff||_{C (\mathbb{R})} \equiv \sup_{\tau \in \mathbb{R}} | (Ff)(\tau)| \le \sqrt \pi ||f||_1.\eqno(2.5)$$
Moreover, if the Mellin transform $(2.1)$ of $f$ satisfies the condition
$${f^*(s)\over \Gamma (1-s)} \in L_p (1-\nu-i\infty, 1-\nu+ i\infty),\  0 < \nu < {1\over 2},\  1 < p \le 2,\eqno(2.6)$$
Then 
$$(Ff)(\tau)=  {2 \sqrt\pi\over  \cosh(\pi\tau)} \int_0^\infty  {\rm Re} \left[ I_{i\tau} \left(\sqrt x\right) \right]  K_{i\tau}(\sqrt x) \varphi (x)dx,\eqno(2.7)$$
where
$$\varphi(x)= {1\over 2\pi i}  \int_{1-\nu- i\infty}^{1-\nu+i\infty} {f^*(s)\over \Gamma(1-s)}   x^{-s} ds,\eqno(2.8)$$
 integral $(2.8)$ converges with respect to the norm in $L_{1-\nu, q}(\mathbb{R}_+),\  q=p/(p-1)$ and the index transform $(2.7)$ is the Lebedev transform with the modified Bessel functions as the kernel \cite{square}.}

\begin{proof} The proof of the norm inequality (2.5) is straightforward from  (1.1) and  inequality (1.10).  The continuity of $(Ff ) (\tau)$ follows from the absolute and uniform convergence of the corresponding integral. In fact, we derive
$$ | (Ff)(\tau)| \le {\sqrt\pi\over  \cosh(\pi\tau)} \int_0^\infty e^{-x/2} \left| I_{i\tau} \left({x\over 2}\right) \right|  |  f(x)| dx $$
$$\le \sqrt{ {\pi\over  \cosh(\pi\tau)}}  \  \sqrt{ {\tanh (\pi\tau) \over  \pi\tau}} \int_0^\infty  |  f(x)| dx  \le \sqrt \pi ||f||_1.$$
Hence we arrive at (2.5).  Further, the asymptotic formulae for  the modified Bessel function of the first kind (1.4), (1.5) and the Macdonald function (1.8)  \cite{erd}, Vol. II yield that for each fixed $\tau > 0$ the kernel in (2.7)  has the behavior 
$$ {\rm Re} \left[ I_{i\tau} \left(\sqrt x\right) \right]  K_{i\tau}(\sqrt x) = O(\log x ),\ x \to 0,$$
$$ {\rm Re} \left[ I_{i\tau} \left(\sqrt x\right) \right]  K_{i\tau}(\sqrt x)=  O\left( {1\over \sqrt x}\right),\ x \to \infty.$$
Hence this  kernel belongs to the space $L_{\nu, p}(\mathbb{R}_+),\  0 < \nu < {1\over 2},\  1 < p \le 2.$   Now by condition (2.6) and Theorem 86 in \cite{tit} we have (see (2.8)) $\varphi \in  L_{1-\nu, q}(\mathbb{R}_+),\  q=p/(p-1).$  

In the meantime, recalling  (1.4),  (1.5), we find 
$$e^{-x/2} {\rm Re} \left[ I_{i\tau} \left({x\over 2}\right) \right] = O(1),\ x \to 0,$$
$$e^{-x/2} {\rm Re} \left[ I_{i\tau} \left({x\over 2}\right) \right] = O\left( {1\over \sqrt x}\right),\ x \to \infty.$$
Hence  the kernel of (1.1) belongs to $L_{\nu, p}(\mathbb{R}_+)$ and $f^*(1-s) \in L_p(\nu-i\infty,  \nu+ i\infty)$ via condition (2.6). Indeed,  we have 
$$\int_{\nu-i\infty}^{ \nu  +i\infty} |f^*(1-s) |^p |ds| =  \int_{1-\nu-i\infty}^{1- \nu  +i\infty} |\Gamma (1-s) f^*(s) |^p {|ds|\over |\Gamma(1-s)|^p} $$ 
$$ \le \left[\Gamma(\nu)\right]^p    \int_{1-\nu-i\infty}^{1- \nu  +i\infty} |f^*(s) |^p {|ds|\over |\Gamma(1-s)|^p} < \infty. $$ 
This means, that $f \in L_{1-\nu,q}(\mathbb{R}_+)$ and integral (1.1) converges absolutely.  Moreover, Theorem 87 in \cite{tit}, the Parseval identity  (2.2), integral representation (1.12) and relation (8.4.23.23) in \cite{prud}, Vol. III  lead to the equalities 
$$(Ff)(\tau)= {1\over 2\pi  i} \int_{\nu-i\infty}^{\nu  +i\infty} \frac {\Gamma(s+ i\tau)\Gamma(s-i\tau) \Gamma(1/2-s)}{\Gamma(s)\Gamma (1-s)} f^*(1-s) x^{-s} ds$$
$$= {2 \sqrt\pi\over  \cosh(\pi\tau)} \int_0^\infty  {\rm Re} \left[ I_{i\tau} \left(\sqrt x\right) \right]  K_{i\tau}(\sqrt x) \varphi (x)dx,$$
where $\varphi$ is defined by (2.8) and both integrals converge absolutely.  This gives (2.7) and completes the proof of Theorem 1.   

\end{proof}

The inversion formula for the transform (1.1) is established by 

{\bf Theorem 2.}  {\it Under conditions of Theorem 1 let also the Mellin transform $f^*(s)$ be analytic in the strip $1/2 < {\rm Re}\  s  <  3/ 2$ and 
$${f^*(s)\over \Gamma (1-s)} \in L_p (1-\nu-i\infty, 1-\nu+ i\infty) \cap L_1 (1-\nu-i\infty, 1-\nu+ i\infty),\ |\nu| < {1\over 2}. \eqno(2.9)$$
Then, assuming that the index transform $(1.1)$ satisfies the integrability condition $(Ff)(\tau) \in L_1(\mathbb{R}_+; \tau e^{\pi\tau} d\tau)$, it has the following inversion formula for all $x >0$ 
$$x f(x)=  { 2  \over \pi\sqrt\pi }   \int_0^\infty  \cosh(\pi\tau)  \left[  \  {}_2F_2  \left( 1,\ {1\over 2} ; \  1 + i\tau; \ 1 - i\tau; \   x \right)  - {\pi\tau  \  e^{x/ 2} \over \sinh(\pi\tau)}  \ {\rm Re} \left[ I_{ i\tau} \left( {x\over 2} \right) \right]   \right]  (Ff)(\tau) d\tau,\eqno(2.10)$$
where the corresponding  integral converges absolutely.}

\begin{proof}   The additional condition (2.9) and analyticity of $f^*(s)$ in the strip ${\rm Re}\  s= 1-\nu \in (1/2, 3/2)$ imply via the Cauchy theorem that integral (2.8) does not change when we shift the contour within the strip.  Meanwhile,  the Lebedev expansion theorem \cite{square} for the index transform (2.7) says that if $\varphi \in L_{3/4,1} \left( (0,1)\right) \cap L_{5/4,1}\left( (1,\infty)\right)$, then for all $x > 0$ the following inversion formula holds
$$\int_x^\infty   \varphi(y) dy = {2\over \pi^2\sqrt\pi } \int_0^\infty \tau \sinh(2\pi\tau) K^2_{i\tau} (\sqrt x) (Ff)(\tau) d\tau.\eqno(2.11)$$
However,  since $\varphi(x) \in L_{1-\nu, q}(\mathbb{R}_+)$,  we let $1/4 < \nu < 1/2$ and use  the H\"{o}lder inequality to get the estimate 
$$ \int_0^1  |\varphi (x)|  x^{- 1/4} dx\le ||\varphi||_{1-\nu,q} \left(\int_0^1  x^{(\nu- 1/4)p-1} dx\right)^{1/p} = {||\varphi||_{1-\nu,q}  \over [(\nu-1/4) p]^{1/p} } < \infty,$$
which guarantees the assumption  $\varphi  \in L_{3/4,1} \left( (0,1)\right)$.   On the other hand,  letting $-1/2 < \nu < -1/4$, we find 
$$\int_1^\infty   |\varphi (x)|  x^{1/4} dx\le ||\varphi||_{1-\nu,q} \left(\int_1^\infty   x^{(\nu+ 1/4)p-1} dx\right)^{1/p} = {||\varphi||_{1-\nu,q}  \over [- (\nu+1/4) p]^{1/p} } < \infty$$
and we have  $\varphi \in L_{5/4,1}\left( (1,\infty)\right).$ Therefore, substituting the value of $\varphi$ by  the integral (2.8) in the left-hand side of (2.11), we change the order of integration,  appealing to  Fubini's theorem. Hence, calculating the elementary inner integral, using the reduction formula for the gamma-function and elementary substitutions, we deduce
$$\int_x^\infty   \varphi(y) dy = - {1\over 2\pi i}  \int_{1-\nu- i\infty}^{1-\nu+i\infty} {f^*(s)\over \Gamma(2-s)}   x^{1-s} ds ={1\over 2\pi i}  \int_{\nu- i\infty}^{\nu+i\infty} {f^*(1-s)\over \Gamma(1+s)}   x^{s} ds,\ -{1\over 2} < \nu < 0.$$
Therefore, combining with (2.11), we arrive at the equality  
$$ {1\over 2\pi i}  \int_{\nu- i\infty}^{\nu+i\infty} {f^*(1-s)\over \Gamma(1+s)}   x^{s} ds =
{ 2\over \pi^2\sqrt\pi } \int_0^\infty \tau \sinh(2\pi\tau) K^2_{i\tau} (\sqrt x) (Ff)(\tau) d\tau.\eqno(2.12)$$
The integral in the right-hand side of (2.12) converges absolutely due to the imposed condition $(Ff)(\tau) \in L_1(\mathbb{R}_+; \tau e^{\pi\tau} d\tau)$ and the Lebedev inequality for the Macdonald function (cf. \cite{yal}, p.  99)
$$|K_{i\tau} (x)| \le {x^{-1/4}\over \sqrt {\sinh(\pi\tau)}},\ x, \tau > 0.\eqno(2.13)$$
Indeed, we have for all $x >0$
$$\int_0^\infty \tau \sinh(2\pi\tau) K^2_{i\tau} (\sqrt x) \left| (Ff)(\tau) \right| d\tau \le  2 x^{-1/4}  \int_0^\infty \tau   e^{\pi\tau} \left| (Ff)(\tau) \right| d\tau < \infty.\eqno(2.14)$$
Hence, returning to (2.12), we apply to the both sides of this equality the Laplace transform with respect to $x$ \cite{tit}.  Changing the order of integration in the left-hand side of the obtained equality  by the Fubini theorem and calculating the inner integral, it gives the result
$${1\over 2\pi i}  \int_0^\infty  e^{-xy}  \int_{\nu- i\infty}^{\nu+i\infty} {f^*(1-s)\over \Gamma(1+s)}   y^{s} ds dy 
= {1\over 2x \pi i}  \int_{\nu- i\infty}^{\nu+i\infty} f^*(1-s) x^{- s} ds=  x^{-2} \ f\left({1\over x}\right),\  x >0 $$ 
by virtue of (2.3) and the condition $f^*(s) \in L_1 (1-\nu-i\infty, 1-\nu+ i\infty)$, which follows immediately from (2.9).  Therefore, for all $x >0$  we derive the formula 
$${1\over x^{2}}  \  f\left({1\over x}\right)  = { 2\over \pi^2\sqrt\pi } \int_0^\infty  e^{-xy}  \int_0^\infty \tau \sinh(2\pi\tau) K^2_{i\tau} (\sqrt y) (Ff)(\tau) d\tau dy$$
$$= { 2\over \pi^2\sqrt\pi }   \int_0^\infty \tau \sinh(2\pi\tau) K(x,\tau)  (Ff)(\tau) d\tau,\eqno(2.15)$$
where
$$K(x,\tau) =  \int_0^\infty  e^{-xy} \  K^2_{i\tau} (\sqrt y) \ dy,\  x, \tau >0\eqno(2.16)$$
and the  interchange of the order of integration is allowed by Fubini's theorem and the estimate (2.14).  Our final goal is to calculate the integral (2.16), which is absent in \cite{prud}, Vol. II.  To do this we appeal again to the Parseval equality (2.2) and relation (8.4.23.27) in \cite{prud}, Vol. III and the Slater residue theorem \cite{prud}, Vol. III.  Hence we deduce $(0 < \gamma < 1)$ 
$$\int_0^\infty  e^{-xy} \  K^2_{i\tau} (\sqrt y) \ dy =  {1\over 4x i\sqrt \pi }  \int_{\gamma- i\infty}^{\gamma +i\infty} 
\Gamma(s+ i\tau)\Gamma(s-i\tau) {\Gamma(s) \Gamma(1-s) \over \Gamma (1/2 +s)}\  x^s ds $$
$$= {\pi\over 2x \sinh(\pi\tau) } \left[ {i\over 2}  \  (4x) ^{-i\tau} \  \Gamma(-i\tau) \  {}_1F_1 
\left( {1\over 2}  + i\tau; \ 1 +2i\tau; \ {1\over x} \right) -  {i \over 2}  \  (4x)^{i\tau} \   \Gamma(i\tau) \right.$$
$$\left. \times \  {}_1F_1  \left( {1\over 2}  - i\tau; \ 1 - 2i\tau; \ {1\over x} \right)  + {1  \over \tau  }
\  {}_2F_2  \left( 1,\ {1\over 2} ; \  1 + i\tau; \ 1 - i\tau; \ {1\over x} \right)   \right]. $$
However, employing relation (7.11.1.5) in \cite{prud}, Vol. III, the latter Kummer functions can be expressed in terms of the modified Bessel functions of the first kind.  Precisely, we obtain
$${}_1F_1 \left( {1\over 2}  \pm  i\tau; \ 1 \pm 2i\tau; \ {1\over x} \right) =  \Gamma(1\pm i\tau) (4x)^{\pm i\tau} e^{1/ (2x)} I_{\pm i\tau} \left( {1\over 2x} \right).\eqno(2.17)$$
Hence after substitution of this expression and straightforward simplifications we get finally the value of the kernel (2.16), namely, 
$$K(x,\tau)=  - \  {\pi\over 2x \sinh(\pi\tau) } \left[ {\pi \  e^{1/ (2x)} \over \sinh(\pi\tau)}  \ {\rm Re} \left[ I_{ i\tau} \left( {1\over 2x} \right) \right]   -  {1  \over \tau  }
\  {}_2F_2  \left( 1,\ {1\over 2} ; \  1 + i\tau; \ 1 - i\tau; \ {1\over x} \right)   \right]. $$
Returning to (2.15) and changing $1/x$ on $x$,  we end up with the inversion formula (2.10), completing the proof of Theorem 2. 
 
\end{proof}

\section{The index transform (1.2)} 

In this section we will examine the boundedness and invertibility conditions for  the operator (1.2). As we see,  it represents a different  transformation, where the integration is realized with respect to the index of the modified Bessel function of the first kind. Such integrals are, generally, unusual and have no common technique to evaluate.  This is why the  mapping properties of the transform (1.2)  and its inversion formula could give such a method of their  evaluation and a source of  new formulas. 

We begin with 

{\bf Theorem 3.}  {\it  The index transform  $(1.2)$  is well-defined as a  bounded operator 
$$G:  L_1(\mathbb{R};\   [\cosh (\pi\tau)]^{-1/2} d\tau ) \to C (\mathbb{R}_+)$$
 and the following norm inequality holds 
$$||Gg ||_{C (\mathbb{R}_+)}  \le \sqrt \pi ||g ||_{L_1(\mathbb{R}; [\cosh (\pi\tau)]^{-1/2} d\tau )}.\eqno(3.1)$$
Besides, if $(Gg ) (x) \in L_{\nu,1} (\mathbb{R}_+),\ 0 < \nu < 1/2$, then for all $x >0$ }
$$\int_0^\infty { (Gg)(t)\over x+t} \ dt = \sqrt \pi \ e^{x/2} \int_{-\infty}^\infty  K_{i\tau} \left({x\over 2} \right)  {g(\tau) \over \cosh(\pi\tau) } d\tau. \eqno(3.2)$$

\begin{proof}  In fact, taking (1.2) and using (1.10), we find

$$ | (Gg)(x)| \le \sqrt\pi \  e^{-x/2} \int_{-\infty}^\infty   \left| I_{i\tau} \left({x\over 2}\right)\right|  \   { |g(\tau)| d\tau\over \cosh(\pi\tau)}\le  \sqrt{ \pi}  \int_{-\infty}^\infty   \sqrt{ {\tanh (\pi\tau) \over  \pi\tau}}  \   { |g(\tau)| d\tau\over \sqrt {\cosh(\pi\tau)}} $$

$$ \le \sqrt \pi ||g ||_{L_1(\mathbb{R}; [\cosh (\pi\tau)]^{-1/2} d\tau )}.$$
This proves (3.1).   When $(Gg ) (x) \in L_{\nu,1} (\mathbb{R}_+),\ 0 < \nu < 1/2$, we take the Mellin transform (2.1) from both sides of (1.2) and change the order integration by Fubini's theorem in the right-hand side of the obtained equality.   Indeed,  from the definition (1.6) of the modified Bessel function of the first kind and similar to (1.7),  (1.10) we derive the inequality 
$$\left|I_{i\tau} \left({x\over 2}\right) \right| \le   I_0  \left({x\over 2}\right)  \sqrt{ {\sinh  (\pi\tau) \over  \pi\tau}},\quad  x >0,\ \tau \in \mathbb{R}.\eqno(3.3)$$
Therefore, taking into account asymptotic formulas (1.4), (1.5), we have  
$$\int_0^\infty |x^{s-1} |  e^{-x/2} \int_{-\infty}^\infty  \left| {\rm Re}  \left[ I_{i\tau} \left({x\over 2}\right)\right] \right| \   { |g(\tau)| d\tau\over \cosh(\pi\tau)} dx$$
$$ \le \int_0^\infty x^{\nu -1}   e^{-x/2}  I_0  \left({x\over 2}\right) dx \int_{-\infty}^\infty  { |g(\tau)| d\tau\over \sqrt{\cosh(\pi\tau)}} < \infty,\  0 < \nu < 1/2.$$
Hence, appealing to Lemma 1, we get 
$$\Gamma(s)\Gamma (1-s) (G g)^* (s)  =   \Gamma(1/2-s) \int_{-\infty}^\infty  \Gamma(s+ i\tau)\Gamma(s-i\tau)  \   g(\tau) d\tau.\eqno(3.4)$$
Meanwhile,  relation (8.4.2.5)  in \cite{prud}, Vol. III and the Parseval identity (2.2), which  still holds for the limit case $p=1$ under conditions of the theorem yield 
$${1\over 2\pi i} \int_{\nu-i\infty}^{\nu+i\infty} \Gamma(s)\Gamma (1-s) (G g)^* (s)  x^{-s} ds = \int_0^\infty { (Gg)(t)\over x+t} \ dt,\  x >0.\eqno(3.5)$$
Consequently, employing relation  (8.4.23.5) in \cite{prud}, Vol. III, we obtain from (3.4), (3.5) 
$$\int_0^\infty { (Gg)(t)\over x+t} \ dt =  {1\over 2\pi i} \int_{\nu-i\infty}^{\nu+i\infty}   \Gamma(1/2-s)  x^{-s} \int_{-\infty}^\infty  \Gamma(s+ i\tau)\Gamma(s-i\tau)  \   g(\tau) d\tau ds $$
$$= \sqrt \pi \ e^{x/2} \int_{-\infty}^\infty  K_{i\tau} \left({x\over 2} \right)  {g(\tau) \over \cosh(\pi\tau) } d\tau,$$
where the interchange of the order of integration in the right-hand side of the latter equality is justified by Fubini's theorem with the use of the inequality (1.100) for the Macdonald function in \cite{yak}, p. 15,  relation (8.4.23.1) in \cite{prud}, Vol. III and  the estimate 
$$\int_{\nu-i\infty}^{\nu+i\infty}  \left| \Gamma(1/2-s)  x^{-s} \right| \int_{-\infty}^\infty  \left| \Gamma(s+ i\tau)\Gamma(s-i\tau)  \   g(\tau) \right| d\tau |ds| $$

$$=   2 \int_{\nu-i\infty}^{\nu+i\infty}  \left| \Gamma(1/2-s)  x^{-s} \right| \int_{-\infty}^\infty   \left| \int_0^\infty K_{2i\tau} (2\sqrt y) y^{s-1} dy \right|   \  |  g(\tau)| d\tau |ds|$$

$$ \le  2 x^{-\nu} \int_{\nu-i\infty}^{\nu+i\infty}  \left| \Gamma(1/2-s)  ds \right|  \int_0^\infty K_{0} \left(2\cos \delta \sqrt y\right) y^{\nu-1} dy  \int_{-\infty}^\infty  e^{-2\delta |\tau|}   \  |  g(\tau)| d\tau < \infty,$$
where $0 < \nu < 1/2$ and $\delta$ is chosen from the interval $[\pi/4,\ \pi/2)$ to satisfy the condition $ g \in L_1(\mathbb{R}; [\cosh (\pi\tau)]^{-1/2} d\tau )$.    Thus we established equality (3.2) and completed the proof of Theorem  3. 

\end{proof} 

The inversion theorem for the index transform (1.2) is given by the following result.

{\bf Theorem 4}.  {\it  Let $g(z/i)$ be an even analytic function in the strip $D= \left\{ z \in \mathbb{C}: \ |{\rm Re} z | < \alpha < 1\right\} ,\  g(0)=g^\prime (0)=0, \  g(z/i)= o(1),\  |{\rm Im} z |\to \infty$ uniformly in $D$.  Then under conditions of Theorem 3  for all  $x \in \mathbb{R}$ the  inversion formula holds for the index transform (1.2)
$$ g(x) =  \lim_{\varepsilon \to 0}  \  {\cosh(\pi x) \over \sqrt \pi }  \int_0^\infty  \left[ {|\Gamma(\varepsilon-1 + ix) |^2  x\sinh (\pi x) \over   \pi \sqrt \pi \ \Gamma(\varepsilon - 1/2)}  \  {}_2F_2  \left(1,  {3\over 2}- \varepsilon; \  2-\varepsilon - ix; \ 2-\varepsilon + ix; \  t \right) \right.$$

$$\left. -   x \   e^{t/2}\  t^{\varepsilon -1}  {\rm Im} \left[ { I_{ix} \left( t/2 \right)  \over \sin(\pi( \varepsilon +ix))} \right] \   \right]   (Gg)(t) \ dt ,\eqno(3.6)$$
where the limit is pointwise.}

\begin{proof}   Indeed,  recalling  (3.2), we multiply its both sides by $ e^{-y/2} K_{ix} \left({y/2} \right) y^{\varepsilon -1}$ for some positive $\varepsilon \in (0,1)$ and integrate with respect to $y$ over $(0, \infty)$.  Hence we obtain 
$$\int_0^\infty  e^{-y/2} K_{ix} \left({y\over 2} \right) y^{\varepsilon -1}  \int_0^\infty { (Gg)(t)\over y +t} \ dt dy  = \sqrt \pi \ \int_0^\infty  K_{ix} \left({y\over 2} \right) y^{\varepsilon -1} \int_{-\infty}^\infty  K_{i\tau} \left({y\over 2} \right)  {g(\tau) \over \cosh(\pi\tau) } d\tau dy. \eqno(3.7)$$
The interchange of the order of integration in the left-hand side of (3.7) can be motivated by Fubini's theorem, employing the generalized Young inequality
$${a^p\over p} + {b^q\over q} \ge  { (ab)^r\over r},\quad   a, b, p, q, r > 0,\ {1\over p} + {1\over q} ={1\over r}.\eqno(3.8)$$
In fact,   from (1.2), (1.10) and the condition $g \in L_1(\mathbb{R};\   [\cosh (\pi\tau)]^{-1/2} d\tau ) $ it follows that $(Gg)(t)= O(1),\ t \to 0$. Now,  letting  in (3.8) $a= (y/q)^{1/p},\  b= (t/p)^{1/q}$, we get 
$$ y+ t \ge \left( \left( {p\over q}\right) ^{q/ (p+q)} +  \left( {q\over p}\right)^{p/(p+q)}\right) y^{ q/ (p+q)} t^ {p/(p+q)}.$$
Hence, choosing positive $p, q$ such that $0< q/(p+q ) < \varepsilon$, we  find the estimate
$$\int_0^\infty  e^{-y/2} \left|K_{ix} \left({y\over 2} \right)\right| y^{\varepsilon -1}  \int_0^\infty { |(Gg)(t) |\over y +t} \ dt dy= \int_0^\infty    e^{-y/2} \left|K_{ix} \left({y\over 2} \right)\right|  y^{\varepsilon -1}  \int_0^1 { |(Gg)(t) |\over y +t} \ dt dy$$

$$+ \int_0^\infty  e^{-y/2}\left| K_{ix} \left({y\over 2} \right) \right| y^{\varepsilon -1}  \int_1^\infty { |(Gg)(t)| \over y +t} \ dt dy
\le C \int_0^\infty    e^{-y/2} K_{0} \left({y\over 2} \right) y^{\varepsilon - q/ (p+q) -1} dy  \int_0^1 t^{-  p/(p+q)} \  dt $$

$$+  \int_0^\infty  e^{-y/2} K_{0} \left({y\over 2} \right) y^{\varepsilon -1} dt  \int_1^\infty  |(Gg)(t) | t^{\nu-1}  dt < \infty,$$
where $C$ is an absolute positive constant.  Thus it guarantees the change of the order of integration in the left-hand side of (3.7).  Then calculating the inner integral, employing formula (2.16.7.6) in \cite{prud}, Vol. II and relation (2.17),  we derive the equality 
$$ \int_0^\infty  K_{ix} \left({y\over 2} \right) y^{\varepsilon -1} \int_{-\infty}^\infty  K_{i\tau} \left({y\over 2} \right)  {g(\tau) \over \cosh(\pi\tau) } d\tau dy =  \int_0^\infty  \left[ {|\Gamma(\varepsilon-1 + ix) |^2\over \Gamma(\varepsilon - 1/2)} 
\  {}_2F_2  \left(1,  {3\over 2}- \varepsilon; \  2-\varepsilon - ix; \ 2-\varepsilon + ix; \  t \right) \right.$$

$$\left.  + { \pi \sqrt \pi  \   e^{t/2}\  t^{\varepsilon -1}\over 2i  \sinh(\pi x)  } \left( { I_{-ix} \left( t/2 \right)  \over \sin(\pi( \varepsilon -ix))} \  -   { I_{ix} \left( t /2 \right)\over  \sin(\pi(\varepsilon+ ix)) }  \right) \right]   (Gg)(t) \ dt.\eqno(3.9)$$
In the meantime,  employing  (1.9) and  the evenness of $g$, the left-hand side of (3.9) can be written as follows
$$ \int_0^\infty  K_{ix} \left({y\over 2} \right) y^{\varepsilon -1} \int_{-\infty}^\infty  K_{i\tau} \left({y\over 2} \right)  {g(\tau) \over \cosh(\pi\tau) } d\tau dy$$

$$= 2 \pi i \int_0^\infty  K_{ix} \left({y\over 2} \right) y^{\varepsilon -1} \int_{-i\infty}^{i\infty}   I_{ z} \left({y\over 2} \right)  {g(z/i) \over \sin (2\pi z) } dz\  dy. \eqno(3.10)$$
On the other hand, according to our assumption $g(z/i)$ is analytic in the vertical  strip $0\le  {\rm Re}  z < \alpha$,  $g(0)=g^\prime (0)=0$ and tends  to zero when $|{\rm Im} | \to \infty $ uniformly in the strip.  Hence,  appealing to the inequality for the modified Bessel   function of the first  kind  (see \cite{yal}, p. 93)
 $$|I_z(y)| \le I_{  {\rm Re} z} (y) \  e^{\pi |{\rm Im} z|/2},\   0< {\rm Re} z < \alpha,$$
one can move the contour to the right in the latter integral in (3.10). Then 

$$2 \pi i  \ \int_0^\infty  K_{ix} \left({y\over 2} \right) y^{\varepsilon -1} \int_{-i\infty}^{i\infty}   I_{ z} \left({y\over 2} \right)  {g(z/i) \over \sin (2\pi z) } dz\  dy$$

$$= 2 \pi i  \  \int_0^\infty  K_{ix} \left({y\over 2} \right) y^{\varepsilon -1} \int_{\alpha -i\infty}^{\alpha + i\infty}   I_{ z} \left({y\over 2} \right)  {g(z/i) \over \sin (2\pi z) } dz\  dy.\eqno(3.11)$$
Now ${\rm Re} z >0$,  and in the right-hand side of (3.11) it is possible to pass to the limit under the integral sign when $\varepsilon \to 0$ and to change the order of integration due to the absolute and uniform convergence.  Therefore the value of the integral (see relation (2.16.28.3) in \cite{prud}, Vol. II
$$\int_0^\infty K_{ix}(y) I_z(y) {dy\over y} = {1\over x^2 + z^2} $$ 
leads us to the equalities 

$$\lim_{\varepsilon \to 0}  2 \pi i  \ \int_0^\infty  K_{ix} \left({y\over 2} \right) y^{\varepsilon -1} \int_{-i\infty}^{i\infty}   I_{ z} \left({y\over 2} \right)  {g(z/i) \over \sin (2\pi z) } dz\  dy$$

$$=    2 \pi i \ \int_{\alpha -i\infty}^{\alpha + i\infty}   {g(z/i) \over (x^2+ z^2) \sin (2\pi z) } dz =  \pi i 
 \left( \int_{-\alpha +i\infty}^{- \alpha- i\infty}   +   \int_{\alpha -i\infty}^{ \alpha+  i\infty}   \right)  {  g(z/i) \  dz \over (z-ix) \  z \sin(2\pi z)}. \eqno(3.12)$$
Hence conditions of the theorem allow to apply the Cauchy formula in the right-hand side of the latter equality in (3.12).  Thus 
$$\lim_{\varepsilon \to 0}  2 \pi i  \ \int_0^\infty  K_{ix} \left({y\over 2} \right) y^{\varepsilon -1} \int_{-i\infty}^{i\infty}   I_{ z} \left({y\over 2} \right)  {g(z/i) \over \sin (2\pi z) } dz\  dy =  { 2\pi^2 \  g(x) \over  x\sinh (2\pi x)} ,\quad x \in \mathbb{R}.$$
Finally, combining with  (3.9),  we arrived at the inversion formula (3.6) and complete to proof of  Theorem 4. 
 
 \end{proof}
 
 {\bf Remark 1}.  When  the passage to the limit under the integral sign is allowed in (3.6), the inversion formula takes the form ( $(Gg)(t) \equiv G(t)$ )
 $$  g(x) =  { \cosh(\pi x) \over \sqrt \pi }  \int_0^\infty  \left[ { x  \   e^{t/2} \over t \sinh(\pi x) }  {\rm Re}
  \left[  I_{ix} \left( {t\over 2} \right) \right] -   {1 \over  2 \pi (1+ x^2) }  \  {}_2F_2  \left(1,  {3\over 2}; \  2 - ix; \  2 + ix; \  t \right) \right] $$
  $$\times   G(t) \  dt.\eqno(3.13)$$

\section{Initial   value problem}

The Lebedev- Skalskaya type transform  (1.2) can be successfully applied to solve an initial value  problem  for the following second  order partial differential  equation, involving the Laplacian
$$\sqrt{ x^2+ y^2}\   \Delta u + x {\partial u \over \partial x}  + y {\partial  u \over \partial y}  +  {1\over 2}  u  =0, \  (x,y) \in 
\mathbb{R}^2 \backslash \{0\},\eqno(4.1)$$ 
where $\Delta = {\partial^2 \over \partial x^2} +  {\partial^2 \over \partial y ^2}$ is the Laplacian in $\mathbb{R}^2$.   Indeed, writing equation (4.1) in polar coordinates $(r,\theta)$, we end up with the equation  
$$ { \partial^2 u \over \partial r^2}  + \left[ {1\over r} +1 \right]   { \partial u \over \partial r} + {1\over r^2} \ { \partial^2 u \over \partial \theta^2} +   {1\over 2r}\ u   = 0.\eqno(4.2)$$

{\bf Lemma 4.} {\it  Let $g(\tau)  \in L_1\left(\mathbb{R}; e^{ (\beta -1/2)  |\tau|} d\tau\right),\  \beta \in (0, 2\pi)$. Then  the function
$$u(r,\theta)=  \sqrt\pi \  e^{-r/2} \int_{-\infty}^\infty   {\rm Re}  \left[ I_{i\tau} \left({r\over 2}\right)\right] \   {  e^{\theta\tau} g(\tau) d\tau\over \cosh(\pi\tau)},\eqno(4.3)$$
 satisfies   the partial  differential  equation $(4.2)$ on the wedge  $(r,\theta): r   >0, \  0\le \theta <  \beta$, vanishing at infinity.}

\begin{proof} In fact, the proof  follows immediately  from the direct substitution (4.3) into (4.2) and the use of  (1.20).  The necessary  differentiation  with respect to $r$ and $\theta$ under the integral sign is allowed via the absolute and uniform convergence, which can be justified, recalling  the inequality (1.10)  and the integrability condition $g \in L_1\left(\mathbb{R}; e^{ (\beta-1/2)  |\tau|} d\tau\right),\  \beta \in (0, 2\pi)$ of the lemma.  Finally,  the condition $ u(r,\theta) \to 0,\ r \to \infty$  is due to the asymptotic formula (1.4) for the modified Bessel function of the first kind. 
\end{proof}

Now,  as a direct consequence of Theorem 4 and Remark 1, we will formulate the initial value problem for equation (4.2) and give its solution.

{\bf Theorem 5.} {\it Let  $g(x)$ be given by formula $(3.13)$, where $G(t)$ satisfies conditions of Theorem 4. 
Then  $u (r,\theta),\   r >0,  \  0\le \theta < \beta$ by formula $(4.3)$  will be a solution  of the initial value problem for the partial differential  equation $(4.2)$ subject to the initial condition}
$$u(r,0) = G (r).$$

\bigskip
\centerline{{\bf Acknowledgments}}
\bigskip

The work was partially supported by CMUP (UID/MAT/00144/2013), which is funded by FCT(Portugal) with national (MEC) and European structural funds through the programs FEDER, underthe partnership agreement PT2020.

\bibliographystyle{amsplain}

\begin{thebibliography}{10}

\bibitem{yak}   Yakubovich S.  Index transforms.   Singapore:  World Scientific Publishing Company; 1996.

\bibitem{erd}    Erd\'elyi A,  Magnus W,   Oberhettinger  F,   Tricomi FG.  Higher transcendental functions. Vols. I,  II. New  York: McGraw-Hill;  1953.

\bibitem{yal}   Yakubovich S,  Luchko Yu.  The hypergeometric approach to integral transforms and convolutions, Mathematics and its applications.  Vol. 287.  Dordrecht:  Kluwer Academic Publishers Group; 1994.

\bibitem{prud}  Prudnikov AP,  Brychkov  YuA,  Marichev OI. Integrals and series:  Vol. I: Elementary functions. New York:  Gordon and Breach;   1986;   Vol. II:  Special functions. New York: Gordon and Breach;  1986;   Vol. III:  More special functions. New York:   Gordon and Breach; 1990.

\bibitem {tit}  Titchmarsh EC.   An introduction to the theory of Fourier integrals.   New York:  Chelsea; 1986.

\bibitem{square} N.N. Lebedev,   On an integral representation of an arbitrary function in terms of squares of Macdonald functions with imaginary index,  {\it Sibirsk. Mat. Zh.},   {\bf 3}  (1962),  213- 222 (in Russian).







\end{thebibliography}

\end{document}